\documentclass[12pt,reqno]{amsart}
\usepackage{amsmath}
\usepackage{amssymb}
\usepackage{amsthm}
\usepackage{mathrsfs}

\usepackage{latexsym}
\usepackage{amssymb}

\newtheorem{theorem}{Theorem}[section]
\newtheorem{lemma}{Lemma}[section]
\newtheorem{corollary}{Corollary}[section]
\newtheorem{remark}{Remark}[section]
\newtheorem{definition}{Definition}[section]
\newtheorem{proposition}{Proposition}[section]
\newtheorem{example}{Example}[section]
\newtheorem{assumption}{Assumption}[section]
\numberwithin{equation}{section}
\newcommand{\bth}{\begin{theorem}}
\newcommand{\ethe}{\end{theorem}}

\newcommand{\bre}{\begin{remark}}
\newcommand{\ere}{\end{remark}}

\newcommand{\ble}{\begin{lemma}}
\newcommand{\ele}{\end{lemma}}
\newcommand{\bde}{\begin{definition}}
\newcommand{\ede}{\end{definition}}

\newcommand{\bco}{\begin{corollary}}
\newcommand{\eco}{\end{corollary}}

\newcommand{\bpr}{\begin{proposition}}
\newcommand{\epr}{\end{proposition}}

\newcommand{\bexer}{\begin{exercise}}
\newcommand{\eexer}{\end{exercise}}
\newcommand{\breh}{\begin{hint}}
\newcommand{\ereh}{\end{hint}}

\newcommand{\halmos}{\hfill \qed}

\newcommand{\bexam}{\begin{example}}
\newcommand{\eexam}{\end{example}}

\newcommand{\pr} {{\bf Proof.}}

\newcommand{\bfi}{\begin{fig}}
\newcommand{\efi}{\end{fig}}

\newcommand{\beao}{\begin{eqnarray*}}
\newcommand{\eeao}{\end{eqnarray*}\noindent}

\newcommand{\beam}{\begin{eqnarray}}
\newcommand{\eeam}{\end{eqnarray}\noindent}
\newcommand{\E}{\mathbf{E}}
\newcommand{\PP}{\mathbf{P}}

\newcommand{\xto}{x\to\infty}

\newcommand{\bM}{\overline{M}}

\newcommand{\bF}{\overline{F}}

\newcommand{\bG}{\overline{G}}
\newcommand{\bV}{\overline{V}}
\newcommand{\bA}{\overline{A}}

\newcommand{\bbr}{{\mathbb R}}

\newcommand{\bbb}{{\mathbb B}}

\newcommand{\bbn}{{\mathbb N}}

\begin{document}
\title[Infinite-time ruin probability with Brownian perturbations]{Infinite-time ruin probability of a multivariate renewal risk model with Brownian perturbations}

\author[D. G. Konstantinides]{Dimitrios G. Konstantinides}

\address{University of the Aegean, Karlovassi, GR-83 200 Samos, Greece}
\email{konstant@aegean.gr.}

\date{{\small \today}}

\begin{abstract}
We consider a multivariate risk model with common renewal process among the lines of business, and Brownian perturbations. Assuming that the integrated tail distribution of claims is multivariate subexponential, we establish an asymptotic relation for the infinite-time ruin probability, which considered on renewal epochs. A more explicit expression is given in case of claim distribution from multivariate regular variation. The results indicate the insensitivity of the asymptotic behavior of the ruin probability with respect to Brownian perturbations. Furthermore, we show that a multivariate distribution with finite expectation, that belongs to the class of multivariate dominatedly varying and long-tailed distributions, possesses integrated tail distribution from the class of multivariate subexponential distributions, which makes easier the checking of conditions in the theorem. 
\end{abstract}

\maketitle
\vspace{3mm}
\vspace{3mm}

\vspace{3mm}

\section{Introduction}

The problem of insurer's ruin in presence of heavy tailed distributions for the claim sizes is well-known already from the last decade of previous century, see in \cite{embrechts:klueppelberg:mikosch:1997}, \cite{asmussen:albrecher:2010}, \cite{konstantinides:2018} for some monographs on this topic. Given that the modern insurance companies are forced to keep more than one risky portfolios, to face the ruin event, as also the dependence structures among the portfolios, the attention of an increasing number of researchers is focused on multivariate risk models with heavy-tailed distribution of claim sizes, see for example \cite{hult:lindskog:2006}, \cite{gao:yang:2014}, \cite{jiang:wang:chen:xu:2015}, \cite{konstantinides:li:2016}, \cite{li:2016}, \cite{samorodnitsky:sun:2016}, \cite{konstantinides:passalidis:2025c}, \cite{passalidis:2025}, \cite{konstantinides:liu:passalidis:2025} among others.  

In this paper, we consider an insurer who operates $d$-lines of business, with $d \in \bbn$, and these lines share a common counting process. 
Let us assume that the claim vectors $\{ {\bf X}^{(i)}\,,\;i \in \bbn \}$ follow distributions, whose support is included on the nonnegative half-axis, 
and each claim vector ${\bf X}^{(i)}=\left(X_1^{(i)}\,,\ldots,\,X_d^{(i)}\right)^{\top}$ may has zero components, but not all of them equal to zero. 
Let us assume that the claim vectors $\{ {\bf X}^{(i)}\,,\;i \in \bbn \}$ arrive at the moments $\{\tau_i\,,\;i \in \bbn \}$, with $\tau_0=0$, that constitute a counting process $\{N(t)\,,\;t\geq 0\}$, defined as follows
\beao
N(t) = \sup\{ n\in \bbn\;:\;\tau_n\leq t \}\,,
\eeao 
for any $t \geq 0$.

Hereafter, we suppose that $\{N(t)\,,\;t\geq 0\}$ is a renewal process, namely the sequence of the inter-arrival times 
$\{ \theta_i = \tau_i - \tau_{i-1}\,,\;i \in \bbn \}$ is sequence of independent and identically distributed (i.i.d.) nonnegative 
random variables. 

Further, we suppose that the insurer charges premiums from the $d$ lines of business, whose  rate of payments is described by the deterministic vector ${\bf p}=(p_1\,,\ldots,\,p_d)^{\top}$, with $p_i \in (0,\,\infty)$, for any $i=1,\,\ldots,\,d$, while he keeps 
initial capital $x>0$, that is allocated on the $d$ lines of business according to ${\bf b}=(b_1\,,\ldots,\,b_d)^{\top}$, with $b_i > 0$, for any $i=1,\,\ldots,\,d$, and $b_1+\cdots+b_d =1$. 

Finally we add one more source of randomness, which stems either from premiums or from claims, described by a multivariate Brownian motion $\{{\bf B}(t)=\left(B_1(t)\,,\ldots,\,B_d(t)\,\right)^{\top}\,,\;t\geq 0\}$ that has arbitrarily correlated components. Hence, if $\vec{\delta} \geq {\bf 0}$ is a fixed vector, named diffusion coefficient, the insurer's surplus process at moment $t>0$ can be described by the equation
\beam \label{eq.KK.1.1} \notag
&&{\bf U}(t) \\[2mm] \notag
&&:=\left( 
\begin{array}{c}
U_{1}(t) \\ 
\vdots \\ 
U_{d}(t) 
\end{array} \right)=x\,\left( 
\begin{array}{c}
b_{1} \\ 
\vdots \\ 
b_{d}
\end{array} \right) + t\,\left( 
\begin{array}{c}
p_{1} \\ 
\vdots \\ 
p_{d}
\end{array} \right) - \left( 
\begin{array}{c}
\sum_{i=1}^{N(t)} X_{1}^{(i)} \\ 
\vdots \\ 
\sum_{i=1}^{N(t)} X_{d}^{(i)}
\end{array} 
\right) + \left( 
\begin{array}{c}
\delta_{1} \\ 
\vdots \\ 
\delta_{d}
\end{array} \right) \odot \left( 
\begin{array}{c}
B_{1}(t) \\ 
\vdots \\ 
B_{d}(t) 
\end{array} \right)\\[2mm]
&&=x\,{\bf b} + t\,{\bf p} - 
\sum_{i=1}^{N(t)} {\bf X}^{(i)}  + \vec{\delta} \odot {\bf B}(t) \,,
\eeam
where by $\odot$ is denoted the Hadamard product. In risk model \eqref{eq.KK.1.1} we make the following assumption.  

\begin{assumption} \label{ass.KK.1.1}
The sequence $\{{\bf X}^{(i)}\,,\;i\in \bbn\}$ contains i.i.d. random vectors with common distribution $F$, and they are independent of the sequence of arrival moments $\{\tau_i\,,\;i \in \bbn \}$. Further we assume that $\{{\bf B}(t)\,,\;t\geq 0\}$ is independent from all other sources of randomness and has nonnegative expectation, in the sense $\E[B_j(t)] \geq 0$, for any $j=1,\,\ldots,\,d$ and any $t\geq 0$.
\end{assumption}

Assumption \ref{ass.KK.1.1} excludes the time-dependent risk models, which are practical enough in modern actuarial practice, see for example \cite{li:tang:wu:2010}, \cite{li:2016}, \cite{yuan:lu:2023}.

In \cite{samorodnitsky:sun:2016} was examined the ruin probability over infinite time horizon in the risk model \eqref{eq.KK.1.1} with $\vec{\delta} = {\bf 0}$, under the conditions of Assumption \ref{ass.KK.1.1}, for the case, when the multivariate integrated tail distribution belongs to class $\mathcal{S}_A$ of multivariate subexponential distributions, see in Section 2 for definitions. Their result generalizes the paper \cite{embrechts:veraverbeke:1982} with respect to two directions, namely with respect to dimension and with respect to renewal counting process instead of a Poisson one.

In this paper we permit the case $\vec{\delta} > {\bf 0}$ as well, in the risk model 
 \eqref{eq.KK.1.1}, while our findings show the 'insensitivity' of the asymptotic 
behavior of the infinite-time ruin probability with respect to Brownian motion, which is observed frequently when the claim sizes have distribution with heavy tail, see in \cite{veraverbeke:1993}, \cite{chen:wang:wang:2013}, \cite{xu:shen:wang:2025}, \cite{li:2017}, \cite{chen:wang:cheng:yan:2023} among others, 
in uni-variate and bi-variate risk models for finite and infinite time cases. Our result is new even in uni-variate case, since it generalizes \cite[Theorem 1]{veraverbeke:1993}, in some sense. This study is restricted to ruin probability only on the arrival moments of the claims (renewal epochs).

The rest of the paper is organized as follows. In Section 2, after some necessary preliminary results about the multivariate heavy-tailed distribution classes, we present a short overview of the results by \cite{samorodnitsky:sun:2016}, that is needed not only for the understanding of the extensions provided by the main result, but also for its proof. In Section 3, we formulate the main result, together with its proof. In Section 4, we show that if a distribution $F$ belongs to the class of multivariate dominatedly varying and long-tailed distributions, symbolically $(\mathcal{D} \cap \mathcal{L})_\mathscr{R} $, and its expectation is finite, then it holds $F_I \in \mathcal{S}_\mathscr{R}$. We also provide an example of building distribution, that belongs to the class $(\mathcal{D} \cap \mathcal{L})_\mathscr{R} $.

\section{Preliminaries} \label{sec.KK.2}

\subsection{Multivariate heavy-tailed distributions} \label{sec.KK.2.1}

Before introduction of the multivariate, heavy-tailed distribution classes, we need some notational conventions. At first, all the limit relations hold as $\xto$. For two real numbers $a,\,b$, we denote by $a\wedge b := \min\{a,\,b\}$. All the vectors are denoted either by Latin bold script or arrowed Greek script and their dimension is equal to $d \in \bbn$. For two vectors ${\bf x}=(x_1,\,\ldots,\,x_d)^{\top}$, ${\bf y}=(y_1,\,\ldots,\,y_d)^{\top}$, and a positive quantity $\kappa > 0$, we denote as usual ${\bf x}\pm {\bf y}=(x_1\pm y_1,\,\ldots,\,x_d\pm y_d)^{\top}$ and $\kappa\,x=(\kappa\,x_1,\,\ldots,\,\kappa\,x_d)^{\top}$, with ${\bf x}^{\top}$ we denote the inverse of ${\bf x}$, while with ${\bf 0}$ we denote the origin of the axes. For any set $\bbb$ from the space $\bbr^d:=(-\infty,\,\infty)^d$, we denote by $\overline{\bbb}$ its closed hull, by $\partial \bbb$ its border and by $\bbb^c$ its complement. We notice that the set $\bbb$ is called increasing, if for any ${\bf x} \in \bbb$ and ${\bf y} \in \bbr_+^d :=[0,\,\infty)^d$, it holds ${\bf x}+{\bf y} \in \bbb$. For two positive $d$-variate functions ${\bf f}$, ${\bf g}$, and some set $\bbb \in \bbr_+^d \setminus \{{\bf 0} \}$, we denote ${\bf f}(x\,\bbb) \sim c\,{\bf g}(x\,\bbb)$, for some constant $c>0$, if it holds
\beao
\lim \dfrac{{\bf f}(x\,\bbb)}{{\bf g}(x\,\bbb)}=c\,.
\eeao  

For a uni-variate distribution $V$, we denote by $\bV$ its tail, namely it holds $\bV(x) = 1- V(x)$, for any $x \in \bbr$. For any $n \in \bbn$, we denote by $V^{n*}$ the $n$-fold convolution of $V$ with itself.

Let us consider now the multivariate, heavy-tailed distribution classes. Initially, we give these classes only for random vectors with distribution $F$, whose support is on $\bbr_+^d$. For the definition of the multivariate subexponentiality, in \cite{samorodnitsky:sun:2016} was used the following fundamental set family:
\beao
\mathscr{R} = \left\{ A\subsetneq \bbr^d\;:\; A\;\text{ open\,,\; increasing}\,,\;A^c\;\text{convex}\,,\;{\bf 0} \notin \bA \right\}\,.
\eeao  
Let us observe that the $\mathscr{R} $ represents a cone with respect to multiplication with positive constants, namely if $A \in \mathscr{R} $ and $k>0$, then it holds $k\,A \in \mathscr{R} $. Thus, for any $A \in \mathscr{R}$, was proved that if ${\bf X}$ is a random vector with distribution $F$, then 
the random variable 
\beao
Y_A:= \sup \{u\;:\;{\bf X} \in u\,A\}
\eeao
has proper distribution $F_A$, whose tail is given by relation
\beao
\bF_A(x) =\PP({\bf X} \in x\,A) = \PP\left( \sup_{{\bf p}\in I_A} {\bf p}^{\top}\,{\bf X} >x\right)\,,\;x>0
\eeao
for some index set $I_A \subsetneq \bbr^d$, see in \cite[Lemmas 4.3(c), 4.5]{samorodnitsky:sun:2016} for more details. 

Therefore, for some fixed set $A \in \mathscr{R} $, we say that $F$ belongs to class of multivariate subexponential distributions on $A$, symbolically $F \in \mathcal{S}_A$, if $F_A \in \mathcal{S}$, that means for any (or, equivalently, for some) integer $n\geq 2$ it holds
\beao
\lim \dfrac{\overline{F^{n*}_A}(x)}{\bF_A(x)} = n\,.
\eeao 
Along the same line, in \cite{konstantinides:passalidis:2024g} was introduced the distribution class $\mathcal{L}_A$, of the multivariate distribution with long tail on $A$, if $F_A \in \mathcal{L}$, namely for any (or, equivalently, for some) $y>0$ it holds
\beao 
\lim \dfrac{\bF_A(x-y)}{\bF_A(x)} =1\,.
\eeao
We denote by $\mathcal{B}_{\mathscr{R}}:= \bigcap_{A \in \mathscr{R}} \mathcal{B}_A$, with $\mathcal{B} \in \{\mathcal{S},\,\mathcal{L}\}$ the classes of multivariate subexponential and multivariate long-tailed distributions, respectively. In the following proposition we reformulate some known results for the classes $\mathcal{S}_A$ and $\mathcal{L}_A$. 

\bpr \label{prop.KK.2.1}
Let $A \in \mathscr{R}$ be some fixed set, $F$ and $G$ be distributions on $\bbr_+^d$, and 
$\mathcal{B} \in \{ \mathcal{S},\,\mathcal{L}\}$. Then we have
\begin{enumerate}
\item[(i)]
$\mathcal{S}_A \subsetneq \mathcal{L}_A$.
\item[(ii)] 
If $F \in \mathcal{B}_A$, and for some $c \in (0,\,\infty)$ it holds $F(x\,A) \sim c\,G(x\,A)$, then 
$G \in \mathcal{B}_A$.
\item[(iii)]
If $F \in \mathcal{L}_A$, then for any ${\bf a} \in \bbr^d$ it holds
$F(x\,A + {\bf a}) \sim F(x\,A)$.
\end{enumerate}
\epr

Assertion $(i)$ of this proposition is direct consequence of the definitions of 
$\mathcal{S}_A $ and $\mathcal{L}_A$, and the fact that $\mathcal{S} \subsetneq \mathcal{L}$, see for example in \cite[Lemma 2]{chistyakov:1964}. Further, Assertion $(ii)$ is implied by \cite[Proposition 4.12(a)]{samorodnitsky:sun:2016}, for the class $\mathcal{S}_A$ and by \cite[Proposition 2.2(iii)]{konstantinides:passalidis:2024g}. Finally, Assertion $(iii)$ can be found in \cite[Proposition 2.3]{konstantinides:passalidis:2024g}.

Finally we present the class of multivariate, regularly varying distributions, symbolically $MRV$, in its standard version. We recall that for an one-dimensional distribution $V$, we say that it belongs to class of regularly varying distributions, with variation index $\alpha \in (0,\,\infty)$,  symbolically $V \in \mathcal{R}_{-\alpha}$, if 
\beao
\lim \dfrac{\bV(t\,x)}{\bV(x)} = t^{-\alpha}
\eeao
for any $t>0$. For a random vector ${\bf X}$ with distribution $F$, we say that 
it belongs to (standard) $MRV$, if there exists an one-dimensional distribution $V \in \mathcal{R}_{-\alpha}$, and some Radon measure $\mu$, non-degenerate to zero, such that holds the limit
\beam \label{eq.KK.2.3}
\lim \dfrac 1{\bV(x)} \PP({\bf X} \in x\,\bbb) = \mu(\bbb)\,,
\eeam
for any Borel set $\bbb \in [0,\,\infty]^d \setminus \{{\bf 0}\} $, which is such that 
$\mu(\partial \bbb) = 0$. In this case, we denote $F \in MRV(\alpha,\,\mu)$. So, if the class $MRV$ represents the union of all $MRV(\alpha,\,\mu)$ from \cite[Proposition 4.14]{samorodnitsky:sun:2016} in combination with Proposition \ref{prop.KK.2.1}(i) we find that
\beam \label{eq.KK.2.4}
MRV \subsetneq \mathcal{S}_{\mathscr{R}} \subsetneq \mathcal{L}_{\mathscr{R}}\,,
\eeam
and relation \eqref{eq.KK.2.4} remains intact for the $\mathcal{B}_A$, with $\mathcal{B} \in \{\mathcal{S},\,\mathcal{L}\}$, for any $A \in \mathscr{R}$.

The $MRV$ distributions, especially in the standard form, are the most popular 
multivariate heavy-tailed distributions and they have found a wide spectrum of applications in probability theory, see in \cite{resnick:2007}, \cite{buraczewski:damek:mikosch:2016}, \cite{samorodnitsky:2016} among others. However, the regular variation of the marginal distributions, in presence of $MRV$ at claim vectors, restricts the actuarial applications and eventually leads to overpricing of the premiums, when the portfolios follow distribution with moderate heavy tails, which are the most common in insurance practice. 

In \cite{konstantinides:liu:passalidis:2025} can be found a variety of general, distribution examples, that belong in class $\mathcal{S}_A$ but not in class $MRV$. Thus, the study of risk models with claim vectors distributions (or integrated tail distributions) from class $\mathcal{S}_A$, does not offer only mathematical generalization, due to \eqref{eq.KK.2.4}, but have impact to the actuarial practice, since it contains a significant multitude of distributions, that do not belong to $MRV$.

\subsection{Ruin probabilities} \label{sec.KK.2.2}

In multivariate risk models, the ruin probability can be defined by several ways, 
since the sets in which can enter the insurer's surplus and to provoke economic 
instabilities for the company, are of several kinds. We refer the reader to \cite{cheng:yu:2019} and \cite{lu:li:yuan:shen:2024} for more discussions about the different definitions of the ruin probabilities. Here, we use the following assumption for the ruin sets, introduced in \cite{samorodnitsky:sun:2016}.

\begin{assumption} \label{ass.KK.2.1}
Let $L$ be some set, that is open, with $L^c$ convex, $-L$ increasing, ${\bf 0} \in \partial L$ and it holds $x\,L =L$, for any $x > 0$. 
\end{assumption}

\bre \label{rem.KK.2.1}
For any set $L$, as described in Assumption \ref{ass.KK.2.1}, it holds
\beao
A:=({\bf b} -L) \in \mathscr{R}\,,
\eeao
which  indicates the immediate relation of the multivariate distribution classes with the ruin probability. The following two ruin sets are very useful in actuarial practice
\beao
L_1:=\left\{ {\bf y}\;:\; \sum_{i=1}^d y_i < 0 \right\}\,, \qquad L_2:=\left\{ {\bf y}\;:\;  y_i < 0\,,\;\exists \;i=1,\,\ldots,\,d\; \right\}\,.
\eeao
\ere

Hereafter, the set $L$ satisfies Assumption \ref{ass.KK.2.1}. 
We denote by ${\bf U}_0(t)$ the surplus in \eqref{eq.KK.1.1} for $\vec{\delta}={\bf 0}$, only in the present subsection. In \cite{samorodnitsky:sun:2016}, was examined the ruin probability for ${\bf U}_0(t)$, in a renewal risk model, under Assumption \ref{ass.KK.2.1}, and in this case the ruin probability has the form
\beao
\psi_{{\bf b},\,L}^o(x)&:=&\PP({\bf U}_0(t) \in L\,,\;\exists\;t>0)= \PP\left(\sum_{i=1}^{N(t)} {\bf X}^{(i)} -t\,{\bf p} \in x\,{\bf b} -L\,,\;\exists\;t>0 \right)\\[2mm]
&=& \PP\left(\sum_{i=1}^{n} \left({\bf X}^{(i)} -\theta_i\,{\bf p}\right) \in x\,A\,,\;\exists\;n \in \bbn \right)
=\PP\left(\sum_{i=1}^{n} {\bf Z}^{(i)} \in x\,A\,,\;\exists\;n \in \bbn \right)\,,
\eeao                                                        
with $\left\{{\bf Z}^{(i)}:= {\bf X}^{(i)} -\theta_i\,{\bf p}\,,\; i \in \bbn \right\}$, which is sequence of\, i.i.d. copies of ${\bf Z}$, since by Assumption \ref{ass.KK.2.1} we have renewal counting process. We find there also the assumption that $\E[{\bf X}] < \overrightarrow{\infty}$, namely $\E\left[X_j \right] < \infty$,  for any $j = 1,\,\ldots,\,d$. To avoid trivial cases, was accepted the assumption of net profit condition
\beao
{\bf c}= -\E[{\bf Z}] > {\bf 0}\,,
\eeao
and for the integrated tail distribution $F_I \in \mathcal{S}_A$, for any $A \in \mathscr{R}$,  where 
\beao
F_I(x\,A) :=\dfrac 1{\Theta} \int_0^{\infty} F(x\,A+v\,{\bf c})\,dv\,,
\eeao
with
\beao
\Theta:= \int_0^{\infty} F([0,\,\infty)^d+v\,{\bf c})\,dv\,,
\eeao
where $\Theta \in (0,\,\infty)$, since ${\bf 0} < \E[{\bf X}] < \overrightarrow {\infty}$.  Hence, under these conditions was found 
\beao
\psi_{{\bf b},\,L}^o(x) \sim \int_0^{\infty} F(x\,A+v\,{\bf c})\,dv=:H(x)\,.
\eeao
 
In this paper, we examine also the case $\vec{\delta}>{\bf 0}$. As we mentioned before, here the study of the ruin probability is restricted only to arrival moments of the claim vectors. In risk model \eqref{eq.KK.1.1} the ruin probability on renewal epochs is reduced to
\beam \label{eq.KK.2.9} 
&&\psi_{{\bf b},\,L}(x) =\PP\left({\bf U}(\tau_i) \in L\,,\;\exists \; i \in \bbn \right)\\[2mm] \notag
&&= \PP\left(\sum_{i=1}^{n} \left({\bf X}^{(i)} -\theta_i {\bf p} - \vec{\delta} \odot {\bf B}^{(i)}\right) \in x A ,\;\exists \; n \in \bbn \right)=\PP\left(\sum_{i=1}^{n} {\bf W}^{(i)} \in x A\,,\;\exists \; n \in \bbn \right),
\eeam
where $\left\{ {\bf B}^{(i)}:= {\bf B}(\tau_{i})- {\bf B}(\tau_{i-1})\,,\; i \in \bbn \right\}$ represents a sequence of\, i.i.d. copies of the vector ${\bf B}:=(B_1,\,\ldots,\,B_d)^{\top}$, since the process $\{{\bf B}(t)\,,\;t\geq 0\}$ is a Brownian motion (with independent increments) and the counting process $\{N(t)\,,\;t\geq 0\}$ is a renewal one, while these two sources of randomness are independent each other. Further, we denote by 
\beao
\left\{ {\bf W}^{(i)}:= {\bf Z}^{(i)} - \vec{\delta} \odot {\bf B}^{(i)}\,,\; i \in \bbn \right\}\,,
\eeao 
a sequence of i.i.d. copies of ${\bf W}$. Due to the fact that the multivariate Brownian motion $\{{\bf B}(t)\,,\;t\geq 0\}$ has nonnegative expectation of each component, we obtain that to avoid the trivial case $\psi_{{\bf b},\,L}(x) =1$, in relation \eqref{eq.KK.2.9}, it should be satisfied a new net profit condition
\beam \label{eq.KK.2.10} 
 {\bf c}^* = -\E[{\bf W}] = {\bf c}+ \vec{\delta} \odot \vec{\mu} >{\bf 0}\,,
\eeam
where $\vec{\mu}=(\mu_1,\,\ldots,\,\mu_d)^{\top}=\left(\E[B_1],\,\ldots,\,\E[B_d]\right)^{\top}$, with ${\bf B}=(B_1,\,\ldots,\,B_d)^{\top}$. From relation \eqref{eq.KK.2.10} we can easily see that if a pair $(\mu_i,\,\delta_i)$, with $i=1,\,\ldots,\,d$ has positive components, then the net profit condition relaxes in comparison with $\vec{\delta} = {\bf 0}$ case, in the sense that it permits the insurer to reduce the premiums at the $i$-line of business.

Finally, we need also the function
\beam \label{eq.KK.2.11} 
F_I^*(x\,A) =\dfrac 1{\Theta^*} \int_0^{\infty} F(x\,A +v\,{\bf c}^* )\,dv\,,
\eeam 
with 
\beao
\Theta^* = \int_0^{\infty} F([0,\,\infty)^d +v\,{\bf c}^* )\,dv \in (0,\,\infty)
\eeao
for the preliminary lemmas of the main result. Let us note that the condition $\E[{\bf X}] < \overrightarrow{\infty}$ remains.

\bre \label{rem.KK.2.2}
The ruin probability \eqref{eq.KK.2.9} , is obviously smaller than
\beam \label{eq.KK.2.a} 
\psi_{{\bf b},\,L}^*(x):=\PP({\bf U}(t) \in L\,,\;\exists\;t>0)\,,
\eeam
that denotes, in some sense, the ruin probability in model \eqref{eq.KK.1.1}, namely includes also the case when the entrance of the surplus into $L$, happens in the time interval $(\tau_{i-1},\,\tau_i)$, due to the Brownian motion. This is a usual simplification technique for risk models with diffusion terms, see for example in \cite{boutsikas:economides:vaggelatou:2024} and references therein. However, in our case, because of the heavy tail of the distribution of claim vectors, there are indications that the $\psi_{{\bf b},\,L}^*(x)$ and $\psi_{{\bf b},\,L}(x)$ are asymptotically equivalent, see in Remark \ref{rem.KK.3.2} below.
\ere

\section{Main result} \label{sec.KK.3}

Here we provide the main result together with its proof.

\subsection{Insensitivity of ruin probability} \label{sec.KK.3.1}

We remind that, in the following theorem the ruin set $L$, satisfies Assumption \ref{ass.KK.2.1}.

\bth \label{th.KK.3.1}
Let $A =({\bf b} -L) \in \mathscr{R}$ be some fixed set. We consider the risk 
model \eqref{eq.KK.1.1}. Under the Assumption \ref{ass.KK.1.1}, if $\E[{\bf X}]< \overrightarrow{\infty}$, 
 $F_I \in \mathcal{S}_A$, and relation \eqref{eq.KK.2.10} is true, then we have
\beam \label{eq.KK.3.1} 
\psi_{{\bf b},\,L}(x) \sim H(x) \,.
\eeam   
\ethe 

\bre \label{rem.KK.3.1}
Theorem \ref{th.KK.3.1} from one hand side generalizes the \cite[Theorem 5.2]{samorodnitsky:sun:2016}, since it permits also the case $\vec{\delta} > {\bf 0}$, and from the other hand side it shows that the asymptotic behavior of the ruin probability remains unchanged, either with or without the presence of the Brownian motion. This last is referred as 'insensitivity of ruin probability with respect to Brownian perturbations'. Furthermore, if instead of condition $F_I \in \mathcal{S}_A$, we require $F_I \in \mathcal{S}_\mathscr{R}$, then relation \eqref{eq.KK.3.1} is true for all the ruin sets $L$, that satisfy Assumption \ref{ass.KK.2.1}. 
\ere

\bre \label{rem.KK.3.2}
In the uni-variate case, where $A^*=(1,\,\infty)$, $b^*=1$ and $L^*=(-\infty,\,0)$ 
relation \eqref{eq.KK.3.1} is reduced into
\beam \label{eq.KK.3.2} 
\psi_{b^*,\,L^*}(x) \sim \int_0^{\infty} \bF(x+v\,c)\,dv =\dfrac 1c \int_x^{\infty} \bF(v)\,dv \,.
\eeam 
Relation \eqref{eq.KK.3.2} represents a partial generalization of \cite[Theorem 1]{veraverbeke:1993}, 
where was considered model \eqref{eq.KK.1.1} where the process $\{N(t)\,,\;t\geq 0\}$ is Poisson, 
with intensity $\lambda>0$, and the Brownian motion has strictly zero expectation. From \eqref{eq.KK.3.2} 
in case of Poisson process and recalling that $c=\E[Z]$, we reach to the formula
\beao
\psi_{b^*,\,L^*}(x) \sim \dfrac 1{\lambda\,p - \E[X]} \int_x^{\infty} \bF(v)\,dv \,,
\eeao
which coincides with that of \cite[Theorem 1]{veraverbeke:1993}. This last relation, shows that in this concrete case the $\psi_{{\bf b},\,L}(x)$ and $\psi_{{\bf b},\,L}^*(x)$ of \eqref{eq.KK.2.9} and \eqref{eq.KK.2.a} are asymptotically equivalent, which means that the ruin event is asymptotically affected only by the renewal epochs and not in the rest time intervals. Hence, we conjecture, that this is valid for all the cases of Theorem \ref{th.KK.3.1}, due to the heavy-tailed distribution of the claims. So, the discretization method for only the renewal epochs by \cite{boutsikas:economides:vaggelatou:2024}, seems to work well here too, although in that paper the distribution of the claims has light tails.  
\ere

Finally we present a corollary, that is implied by Theorem \ref{th.KK.3.1}, if we require that $F \in MRV(\alpha,\,\mu)$, with $\alpha \in (1,\,\infty)$. It is worth to notice that that if $F \in MRV(\alpha,\,\mu)$, with $\alpha \in (1,\,\infty)$, then we have $F_I \in \mathcal{S}_\mathscr{R}$, which follows easily from the line of proof in \cite[Proposition 5.4]{samorodnitsky:sun:2016}. However, we give a more general example in Section 4, which follows by similar methodology. To avoid the trivial case, where ${\bf X}$ has some degenerate components, we suppose that $\mu({\bf x}\;:\;x_i>0) > 0$, for any $i=1,\,\ldots,\,d$. The following result presents a more explicit asymptotic expression in comparison with \eqref{eq.KK.3.1} and generalizes \cite[Proposition 5.4]{samorodnitsky:sun:2016}. The proof of this corollary follows a similar line with that of  \cite[Proposition 5.4]{samorodnitsky:sun:2016} by just making use of Theorem \ref{th.KK.3.1}, therefore we omit the detailed proof.

\bco \label{cor.KK.3.1}
Let $L$ satisfy Assumption \ref{ass.KK.2.1}. We consider risk model \eqref{eq.KK.1.1}, under Assumption \ref{ass.KK.1.1} and the conditions  $F \in MRV(\alpha,\,\mu)$, with $\alpha \in (1,\,\infty)$, and \eqref{eq.KK.2.10} hold. Then we have 
\beao
\psi_{{\bf b},\,L}(x) \sim x\,\bV(x)\, \int_0^{\infty} \mu({\bf b}-L+v\,{\bf c})\,dv \,.
\eeao
\eco

\subsection{Proof of Theorem \ref{th.KK.3.1} } \label{sec.KK.3.2}

Before the proof of the main result, we need two preliminary lemmas. Let us remind 
the definition of $F_I^*$ by relation \eqref{eq.KK.2.11}.

\ble \label{lem.KK.3.1}
Let $A \in \mathscr{R}$ be some fixed set. If $F_I \in \mathcal{L}_A$, then we find
\beam \label{eq.KK.3.4} 
F_I^*(x\,A) \sim \dfrac{\Theta}{\Theta^*}\,F_I(x\,A) \,.
\eeam
\ele

\pr~
From the one hand side, since $\vec{\delta},\,\vec{\mu} \geq {\bf 0}$, and the set $A$ 
is increasing, we obtain
\beam \label{eq.KK.3.5} 
F_I^*(x\,A) &=& \dfrac 1{\Theta^*}\,\int_0^{\infty} \PP({\bf X} \in x\,A + v\,{\bf c} + v\,\vec{\delta} \odot \vec{\mu}  )\,dv\\[2mm] \notag
&\leq& \dfrac{\Theta}{\Theta^*}\,\dfrac 1{\Theta}\,\int_0^{\infty} \PP({\bf X} \in x\,A + v\,{\bf c} )\,dv= \dfrac{\Theta}{\Theta^*}\,F_I(x\,A) \,,
\eeam
that provides the desired upper bound for relation \eqref{eq.KK.3.4}. 

For the lower bound, from \cite[Lemma 4.3(d)]{samorodnitsky:sun:2016} we can find 
some $u=u(v) \in (0,\,x)$, such that it holds
\beam \label{eq.KK.3.6} 
F_I^*(x\,A) &=& \dfrac 1{\Theta^*}\,\int_0^{\infty} \PP({\bf X} \in x\,A + v\,{\bf c} + v\,\vec{\delta} \odot \vec{\mu}  )\,dv\\[2mm] \notag
&\geq&  \dfrac{\Theta}{\Theta^*}\,\dfrac 1{\Theta}\,\int_0^{\infty} \PP({\bf X} \in (x +u)\,A + v\,{\bf c} )\,dv \sim \dfrac{\Theta}{\Theta^*}\,\dfrac 1{\Theta}\,\int_0^{\infty} \PP({\bf X} \in x\,A + v\,{\bf c} )\,dv\,,
\eeam
where at the last step we used the fact that $F_I \in \mathcal{L}_A$, due to dominated convergence theorem. 
~\halmos

\bre \label{rem.KK.3.3}
Relation \eqref{eq.KK.3.4} implies that if $F_I \in \mathcal{B}_A$, with $\mathcal{B} \in \{\mathcal{S},\,\mathcal{L}\}$, then in combination with Proposition \ref{prop.KK.2.1}(ii), we have $F_I^* \in \mathcal{B}_A$. Lemma \ref{lem.KK.3.1} is a crucial step in simplification of the proof of Theorem \ref{th.KK.3.1}. 
\ere

For next lemma we need the following function
\beao
H^*(x) = \int_0^{\infty} F( x\,A + v\,{\bf c}^* )\,dv\,.
\eeao

\ble \label{lem.KK.3.2}
Let $A \in \mathscr{R}$ be some fixed set. If $F_I \in \mathcal{L}_A$, then we have
\beam \label{eq.KK.3.8} 
H^*(x) \sim H(x)\,.
\eeam
\ele

\pr~
The upper bound in relation \eqref{eq.KK.3.8}, follows by similar arguments with that in 
relation \eqref{eq.KK.3.5}.

For the lower bound in relation \eqref{eq.KK.3.8}, we use the similar line 
with that in relation \eqref{eq.KK.3.6} in combination with the fact 
that $H \in \mathcal{L}_A$, which follows by the inclusion $F_I \in \mathcal{L}_A$, 
and the fact that the $F_I$ and $H$ are strongly equivalent, see Proposition  \ref{prop.KK.2.1}(ii).   
~\halmos

We can now present the proof of the main result.

{\bf Proof of Theorem \ref{th.KK.3.1}}~
Following similar line with the proof of \cite[Theorem 5.2]{samorodnitsky:sun:2016}, taking into account Lemma \ref{lem.KK.3.1}, we easily find that
\beam \label{eq.KK.3.9} 
\psi_{{\bf b},\,L}(x) \sim \int_0^{\infty} \PP({\bf W} \in x\,A + v\,{\bf c}^* )\,dv \,.
\eeam
Now, we shall show that
\beam \label{eq.KK.3.10} 
\int_0^{\infty} \PP({\bf W} \in x\,A + v\,{\bf c}^* )\,dv \sim H^*(x)\,.
\eeam
We notice that through similar steps of relation \cite[Eq. (5.7)]{samorodnitsky:sun:2016} we can verify that
\beam \label{eq.KK.3.11} 
I(x):=\int_0^{\infty} \PP({\bf Z} \in x\,A + v\,{\bf c}^* )\,dv \sim H^*(x)\,.
\eeam
From the fact that $F_I \in \mathcal{S}_A$, and through Lemma  \ref{lem.KK.3.2} we find
\beao
F_I(x\,A) \sim \dfrac 1{\Theta}\,H^*(x)\,,
\eeao
therefore, by Proposition \ref{prop.KK.2.1}(ii) follows that $H^* \in \mathcal{S}_A$, and by relation \eqref{eq.KK.3.11} that $I(x) \in \mathcal{S}_A \subsetneq \mathcal{L}_A$. 

Next, we estimate the lower bound of relation \eqref{eq.KK.3.10}. Let 
\beao
\bM := \max_{1\leq j \leq d} |\delta_j\,B_j|
\eeao
be a random variable. Since the set $A$ is increasing, and through \cite[Lemma 4.3(d)]{samorodnitsky:sun:2016}, for some $u \in (0,\,x)$ we obtain
\beam \label{eq.KK.3.12} \notag
\int_0^{\infty} \PP({\bf W} \in x\,A + v\,{\bf c}^* )\,dv &=& \int_0^{\infty} \PP({\bf Z} -\vec{\delta}\odot {\bf B} \in x\,A + v\,{\bf c}^* )\,dv \\[2mm] \notag
&\geq& \int_0^{\infty} \PP({\bf Z} -\bM \in x\,A + v\,{\bf c}^* )\,dv \\[2mm] \notag
&=& \int_0^{\infty} \int_0^{\infty} \PP({\bf Z} \in x\,A + v\,{\bf c}^* +y)\,\PP(\bM \in dy)\,dv \\[2mm] \notag
&\geq& \int_0^{\infty} \int_0^{\infty} \PP({\bf Z} \in (x+u)\,A + v\,{\bf c}^* )\,\PP(\bM \in dy)\,dv \\[2mm]
&\sim& \int_0^{\infty} \PP({\bf Z} \in x\,A + v\,{\bf c}^* )\,dv \sim H^*(x)\,,
\eeam
where at the next-to-last step we used the property of class $\mathcal{L}_A $ for $I(x)$, while at the last step we considered relation \eqref{eq.KK.3.11}. The previous relation renters the desired upper bound for \eqref{eq.KK.3.10}. 

Further, we deal with the lower bound of relation \eqref{eq.KK.3.10}. Let 
\beao
\underline{M} := \left( \min_{1\leq j \leq d} \delta_j\,B_j\right)\bigwedge 0
\eeao
be a non-positive random variable. Following similar symmetric steps with \eqref{eq.KK.3.12} we find
\beao
&&\int_0^{\infty} \PP({\bf W} \in x\,A + v\,{\bf c}^* )\,dv = \int_0^{\infty} \PP({\bf Z} -\vec{\delta}\odot {\bf B} \in x\,A + v\,{\bf c}^* )\,dv \\[2mm]
&& \leq \int_0^{\infty} \PP({\bf Z} - \underline{M} \in x\,A + v\,{\bf c}^* )\,dv = \int_0^{\infty} \int_0^{\infty} \PP({\bf Z}+y \in x\,A + v\,{\bf c}^*)\,\PP(-\underline{M} \in dy)\,dv \leq \\[2mm]
&&\int_0^{\infty} \int_0^{\infty} \PP({\bf Z} \in (x-u)\,A + v\,{\bf c}^*)\,\PP(-\underline{M} \in dy)\,dv \sim \int_0^{\infty} \PP({\bf Z} \in x\,A + v\,{\bf c}^* )\,dv \sim H^*(x).
\eeao
Hence, \eqref{eq.KK.3.10} is true. Combining  \eqref{eq.KK.3.9} with \eqref{eq.KK.3.10} and \eqref{eq.KK.3.8}, we get relation \eqref{eq.KK.3.1}.
~\halmos

\section{Integrated tails in $\mathcal{S}_\mathscr{R}$}  \label{sec.KK.4}

In this section we show that if a distribution $F$ belongs to a concrete distribution 
class, then for the integrated tail distribution the relation $F_I \in \mathcal{S}_\mathscr{R}$ is true. 
So, we obtain a useful tool for the control of the condition in Theorem \ref{th.KK.3.1}.

Let start with the class $\mathcal{D}$, of dominatedly varying distributions. 
An one-dimensional distribution $V$, is such that $V \in \mathcal{D}$, if $\bV(x)>0$, for all $x>0$, and for any (or, equivalently, for some) $b \in (0,\,1)$ it holds
\beao
\limsup \dfrac{\bV(b\,x)}{\bV(x)} < \infty\,.
\eeao
It is well-known that $\mathcal{D} \cap \mathcal{S} \equiv \mathcal{D} \cap \mathcal{L}$, see \cite{goldie:1978}, and further we find
\beam \label{eq.KK.4.1} 
\bigcup_{0 < \alpha < \infty} \mathcal{R}_{-\alpha} \subsetneq \mathcal{D} \cap \mathcal{L}  \subsetneq \mathcal{S}  \subsetneq \mathcal{L}\,,
\eeam
see in \cite[Chapter 2]{leipus:siaulys:konstantinides:2023} for more details. Thus, 
in \cite{konstantinides:passalidis:2024g} were introduced the corresponding multivariate versions of $\mathcal{D}$ and $\mathcal{D} \cap \mathcal{L}$. For a set $A \in \mathscr{R}$, we say that $F \in \mathcal{D}_A$, if $F_A \in \mathcal{D}$, and further that $F\in (\mathcal{D} \cap \mathcal{L})_A$, if $F_A \in \mathcal{D} \cap \mathcal{L}$. Let us denote 
\beao
\mathcal{D}_\mathscr{R} := \bigcap_{A \in \mathscr{R}} \mathcal{D}_A\,,\qquad (\mathcal{D} \cap \mathcal{L})_\mathscr{R} := \bigcap_{A \in \mathscr{R}} (\mathcal{D} \cap \mathcal{L})_A\,.
\eeao
By \cite[Proposition 2.1]{konstantinides:passalidis:2024g} and relations \eqref{eq.KK.2.4}, \eqref{eq.KK.4.1} we have
\beam \label{eq.KK.4.2} 
MRV \subsetneq (\mathcal{D} \cap \mathcal{L})_\mathscr{R} \subsetneq \mathcal{S}_\mathscr{R} \subsetneq \mathcal{L}_\mathscr{R}\,.
\eeam
Following the line of proof of \cite[Proposition 5.4]{samorodnitsky:sun:2016}, we find that 
if $F \in MRV$, with index $\alpha > 1$ (to secure a finite expectation), it implies that 
$F_I \in \mathcal{S}_\mathscr{R}$. 

The following result shows that if $F \in (\mathcal{D} \cap \mathcal{L})_\mathscr{R}$, with 
$\E[{\bf X}]< \overrightarrow{\infty}$. Then we obtain $F_I \in \mathcal{S}_\mathscr{R}$, 
and hence we get weaker sufficient conditions for the property of multivariate subexponentiality 
of the integrated tail distribution.

\bpr \label{pr.KK.4.1}
Let ${\bf X}$ be nonnegative random vector with distribution $F \in (\mathcal{D} \cap \mathcal{L})_\mathscr{R}$ 
and $\E[{\bf X}]< \overrightarrow{\infty}$. Then it holds $F_I \in (\mathcal{D} \cap \mathcal{L})_\mathscr{R} \subsetneq \mathcal{S}_\mathscr{R}$.
\epr

\pr~
Let $A \in \mathscr{R}$. We can write for the integrated tail distribution
\beam \label{eq.KK.4.3}
\bF_{I,A}(x)=F_I(x\,A)= \dfrac 1{\Theta} \int_0^{\infty} F(x\,A +v\,{\bf c})\,dv\,.
\eeam
Hence, it is enough to show that $F_{I,A} \in \mathcal{D} \cap \mathcal{L}$.
 
Let deal first with class $\mathcal{L}$. Let 
\beao
A^*:=A + v\,{\bf c}\,,
\eeao 
from where, since $A \in \mathscr{R}$, it follows $A^* \in \mathscr{R}$. Indeed, by \cite[Lemma 4.3]{samorodnitsky:sun:2016}, we find that there exists a $u \in (0,\,1)$ such that it holds $A^* \subsetneq (1-u)\,A$, with the right member of this inclusion to belong in $\mathscr{R}$, since $\mathscr{R}$ is closed with respect to positive scalar multiplication.  Let us consider some $y >0$, then we obtain
\beam \label{eq.KK.4.4} \notag
1 &\leq& \lim \dfrac{\bF_{I,A}(x-y)}{ \bF_{I,A}(x)} =\lim \dfrac {\int_0^{\infty} \PP\left({\bf X} \in (x-y)\,A + v\,{\bf c}\right)\,dv}{\int_0^{\infty} \PP\left({\bf X} \in x\,A + v\,{\bf c}\right)\,dv} \\[2mm]
&=& \lim \dfrac {x\int_0^{\infty} \PP\left({\bf X} \in (x-y)\,\left(A + \dfrac x{x-y}\,v\,{\bf c}\right)\,\right)\,dv}{x\,\int_0^{\infty} \PP\left({\bf X} \in x\,(A + v\,{\bf c})\,\right)\,dv}\\[2mm] \notag
&\leq& \lim \dfrac {\int_0^{\infty} \PP\left({\bf X} \in (x-y)\,A^*\,\right)\,dv}{\int_0^{\infty} \PP\left({\bf X} \in x\,A^* \right)\,dv} = \lim \dfrac {\int_0^{\infty} \PP\left({\bf X} \in x\,A^*\,\right)\,dv}{\int_0^{\infty} \PP\left({\bf X} \in x\,A^* \right)\,dv}=1\,,
\eeam
where at the third step we made change of variables via $z=x/v$, while at the next-to-last step we used the fact that $F \in \mathcal{L}_\mathscr{R}$, through the dominated convergence theorem, because of $\Theta^*<\infty$. Hence $F_{I,A} \in \mathcal{L}$, thus we obtain $F_I \in \mathcal{L}_A$.

Now, we continue to show that $F_{I,A} \in \mathcal{D}$. Let $b \in (0,\,1)$, then it holds
\beam \label{eq.KK.4.5} \notag
&& \limsup \dfrac{\bF_{I,A}(b\,x)}{ \bF_{I,A}(x)} =\limsup \dfrac {\int_0^{\infty} \PP\left({\bf X} \in b\,x\,A + v\,{\bf c}\right)\,dv}{\int_0^{\infty} \PP\left({\bf X} \in x\,A + v\,{\bf c}\right)\,dv} \\[2mm]
&&= \limsup \dfrac {x\,b\,\int_0^{\infty} \PP\left({\bf X} \in b\,x\,A + b\,x\,v\,{\bf c}\right)\,dv}{x\,\int_0^{\infty} \PP\left({\bf X} \in x\,(A + v\,{\bf c})\,\right)\,dv}\\[2mm] \notag
&&= b\,\limsup \dfrac {\int_0^{\infty} \PP\left({\bf X} \in b\,x\,A^*\right)\,dv}{\int_0^{\infty} \PP\left({\bf X} \in x\,A^* \right)\,dv} \leq b\,C\, \limsup \dfrac {\int_0^{\infty} \PP\left({\bf X} \in x\,A^*\,\right)\,dv}{\int_0^{\infty} \PP\left({\bf X} \in x\,A^* \right)\,dv}=b\,C < \infty\,,
\eeam
where at the second step we made change of variables $z=b\,x/v$ on the numerator, and $z=x/v$ on the denominator, while at the fourth step the constant $C \in [1,\,\infty)$ is implied by the fact that $F \in \mathcal{D}_\mathscr{R}$, due to the dominated convergence theorem. Therefore, from relation \eqref{eq.KK.4.5} we obtain $F_{I,A} \in \mathcal{D}$, that means $F_{I} \in \mathcal{D}_A$. In combination with the previous remarks we have  $F_{I} \in (\mathcal{D} \cap \mathcal{L})_A$, and due to arbitrariness in the choice of $A \in \mathscr{R}$, it follows that $F_{I} \in (\mathcal{D} \cap \mathcal{L})_\mathscr{R}$. 
~\halmos

\bre \label{rem.KK.4.1}
The question about what distributions $F$ are such that the integrated tail distribution $F_I$ are (multivariate) subexponential, has both theoretical and practical interest. In practice it is much easier to check the distribution class of $F$, in comparison with that of $F_I$. In case of uni-variate distribution this question was examined by several papers, see \cite{embrechts:omey:1984} and \cite{klueppelberg:1988} for more discussions. 
\ere

We close this section providing a simple example of distributions that belong to
$(\mathcal{D} \cap \mathcal{L})_\mathscr{R}$, which is inspired by \cite[Example 4.17]{samorodnitsky:sun:2016}, see also \cite[Example 4.6]{konstantinides:liu:passalidis:2025}. Let us notice that in this example $Y_A$ follows distribution $F_A$, as in previous text.

\bexam \label{exam.KK.4.1}
Let $Z$ be a random variable with distribution $G$, whose tail is given as follows:
\beao
\bG(x) = \PP\left( \sum_{j=1}^d X_j >x \;\bigg|\; \dfrac{{\bf X}}{\|{\bf X}\|} \in (\theta_1,\,\ldots,\,\theta_d)^{\top}\, \right)
\eeao
with $ (\theta_1,\,\ldots,\,\theta_d)^{\top} \in \Delta_d$, where $\Delta_d$ represents a unit sphere in $\bbr^d$, and $\|\cdot\|$ is the $L_1$-norm. Let $A \in \mathscr{R}$, be some fixed set and we introduce the random variable
\beao
W:= \inf \left\{ u>0\;:\; u\,\dfrac{{\bf X}}{\|{\bf X}\|} \in A \right\}\,.
\eeao
If $Z,\,W$ are independent, we define the following random variable
\beao
Y_A\stackrel{d}{=} \dfrac Z{W}\,,
\eeao
where $\stackrel{d}{=}$ is the equality in distribution. Let $G \in \mathcal{D} \cap \mathcal{L}$. Then, 
since the $W$ is bounded away from zero, the $W^{-1}$ is bounded from above. So from 
\cite[Theorem 2.2(iii), 3.3(ii)]{cline:samorodnitsky:1994} we obtain 
$F_A \in \mathcal{D} \cap \mathcal{L}$. Thus, for any $A \in \mathscr{R}$, we find 
$F \in (\mathcal{D} \cap \mathcal{L})_A$, which implies $F \in (\mathcal{D} \cap \mathcal{L})_\mathscr{R}$.~\halmos
\eexam

Another example with distribution $F \in (\mathcal{D} \cap \mathcal{L})_\mathscr{R}$ can be found in \cite[Proposition 4.15]{samorodnitsky:sun:2016}.

\section*{Declarations. }

\subsection*{Publishing policy. }
I have read and understood the publishing policy, and submit this manuscript in accordance with this policy.

\subsection*{Competing interests. }
I declare that the author has no competing interests as defined by Springer, or other interests that might be perceived to influence the results and/or discussion reported in this paper.

\subsection*{Third party material. }
All of the material is owned by the author and/or no permissions are required.

\subsection*{Data availability. }
Not applicable

\subsection*{Research funding. }
This research did not receive funding.

\section*{Acknowledge}
I feel the pleasant duty to express my gratitude to C. D. Passalidis, for several discussions on this topic. Further I am grateful to the anonymous referee for several suggestions that improved the paper.

\end{document}